\documentclass[12pt]{article}
\usepackage{a4wide}
\usepackage{amsmath,amssymb,amsthm}
\usepackage[mathscr]{eucal}
\usepackage{graphics}
\usepackage[hypertex]{hyperref}

\author{Ievgen Bondarenko\thanks{The author was partially supported by NSF grants DMS-0308985 and
DMS-0456185}}
\title{\textbf{Dynamics of piecewise linear maps and sets of nonnegative matrices}}

\newcommand{\mset}{\mathscr{K}}

\newtheorem{teor}{Theorem}
\newtheorem{prop}[teor]{Proposition}
\newtheorem{col}[teor]{Corollary}
\newtheorem{lemm}{Lemma}
\theoremstyle{definition}
\newtheorem{defi}{Definition}

\begin{document}
\maketitle

\begin{abstract}
We consider maps $f_{\mset}(v)=\min_{A\in\mset}{Av}$ and
$g_{\mset}(v)=\max_{A\in\mset}{Av}$, where $\mset$ is a finite set
of nonnegative matrices and by ``min'' and ``max'' we mean
component-wise minimum and maximum. We transfer known results about
properties of $g_{\mset}$ to $f_{\mset}$. In particular we show
existence of nonnegative generalized eigenvectors of $f_{\mset}$,
give necessary and sufficient conditions for existence of strictly
positive eigenvector of $f_{\mset}$, study dynamics of $f_{\mset}$
on the positive cone. We show the existence and construct matrices
$A$ and $B$, possibly not in $\mset$, such that $f_{\mset}^n(v)\sim
A^nv$ and $g_{\mset}^n(v)\sim B^nv$ for any strictly positive vector
$v$.
\end{abstract}

\section{Introduction}
The theory of nonnegative matrices has been very well developed
since its appearance in papers of Perron and Frobenius. Vast number
of applications to dynamic programming, probability theory,
numerical analysis, mathematical economics, fractal geometry raise
even greater interest to this field. As a result there are many
strong generalizations of the Perron-Frobenius theory (see
\cite{bapat,berm:nonmatr,gantmacher,gaubert:hommon,goldb,rubinov,
MandlSeneta,metz,rothbl:alg}).


The classical Perron-Frobenius theorem shows that a nonnegative
matrix has a nonnegative eigenvector associated with its spectral
radius, and if the matrix is irreducible then this nonnegative
eigenvector can be chosen strictly positive.
One important generalization of this result was obtained by
U.G.~Rothblum \cite{rothbl:alg}, who studied the structure of the
algebraic eigenspaces of nonnegative matrices and described the
combinatorics that stands behind the index of the spectral radius
and dimensions of the algebraic eigenspaces. Moreover, it was shown
that the algebraic eigenspace of a nonnegative matrix corresponding
to its spectral radius is spanned by a set of nonnegative
generalized eigenvectors with certain strictly positive entries.

Many generalizations of the Perron-Frobenius theory involve
homogeneous monotone functions, especially functions of the form
\[
g_{\mset}(x)=\max_{A\in\mset} Ax,
\]
where $\mset$ is a set of square nonnegative matrices of fixed
dimension and by ``max'' we mean component-wise maximum. Such
functions appear in many well-known problems, such as the theory of
controlled Markov chains, Leontief substitution systems, controlled
branching processes, parallel computations, transportation networks,
etc. The study of maps $g_{\mset}$ was initiated by Richard Bellman.
Using the Brouwer fixed point theorem he proved existence of a
strictly positive eigenvector of the map $g_{\mset}$ in the case
when each matrix in $\mset$ is positive and studied the asymptotic
behavior of iterations $g_{\mset}^n(v)=g_{\mset}(g_{\mset}(\ldots
g_{\mset}(v)\ldots))$ for a nonnegative vector $v$ (see
\cite{bellman:quasi_l_eq} and \cite[chapter XI,
sections~$10-11$]{bellman:dynpr}). These results were generalized to
a set of irreducible matrices by P.~Mandl and
E.~Seneta~\cite{MandlSeneta}.

The most important results for our investigation were obtained by
W.H.M.~Zijm in \cite{zij84generalized}. He showed that there is a
simultaneous block-triangular decomposition of the set of matrices
$\mset$, which was used to give the necessary and sufficient
conditions for the existence of a strictly positive eigenvector of
$g_{\mset}$ and extend the above mentioned result of U.G.~Rothblum
on nonnegative generalized eigenvectors to $g_{\mset}$. Related
results can be found in \cite[Chapter~35]{whittle:opt_over_time}.
Independently, Karel Sladk{\'y} \cite{sladky:bound1,sladky:bound2}
obtained the same block-triangular decomposition and used it to get
bounds on the asymptotic behavior of iterations $g_{\mset}^n(v)$ for
a nonnegative vector $v$. Stronger results about asymptotic behavior
of iterations $g_{\mset}^n(v)$ were obtained
in~\cite{sladky:bound2,sladky:funceq,zijm:asympt} for the case when
some special matrices in $\mset$ are aperiodic.

We consider maps of similar form, but with ``minimum'' instead of
``maximum'':
\[
f_{\mset}(x)=\min_{A\in\mset} Ax.
\]
These maps appear in~\cite{kigami:short} in connection with the
construction of ``self-similar'' metrics on self-similar sets and
finding their Hausdorff dimensions. Also such maps appear in the
study of growth of Schreier graphs of groups generated by finite
automata~\cite{bond:schreier,bondnekr}. These problems demand us to
study spectral properties of maps $f_{\mset}$ and describe
asymptotic behavior of its iterations.

Considering maps $f_{\mset}$ and $g_{\mset}$ we can always suppose
that the set $\mset$ satisfies the product property, i.e. $\mset$ is
constructed by all possible interchanges of corresponding rows
selected from a finite set of square nonnegative matrices (see the
precise definition and explanation in
Section~\ref{section_exist_str_pos_f}). Under this assumption, for
every vector $v$ there exist $A=A_v\in\mset$ and $B=B_v\in\mset$
such that $f_{\mset}(v)=Av$ and $g_{\mset}(v)=Bv$. In the theory of
Markov decision processes this property is usually called the
optimal choice property (see \cite{howard:dynpr,bellman:dynpr}).

The asymptotic behavior of iterations $f^n(v)$ is studied with
respect to the following equivalence. Let $a_n, b_n$, $n\geq 1$, be
sequences of nonnegative numbers or vectors of the same dimension.
We say that $a_n\preceq b_n$ if there exists constant $q>0$ such
that $a_n\leq q\cdot b_n$ for all $n$ large enough. If $a_n\preceq
b_n$ and $b_n\preceq a_n$ then we say that $a_n\sim b_n$ and that
$a_n$ and $b_n$ have the same \emph{growth}. Then $h^{n}(v)\sim
h^{n}(u)$ for any homogenous nondecreasing function
$h:\mathbb{R}^{N}_{+}\rightarrow\mathbb{R}^{N}_{+}$ and any strictly
positive vectors $v,u$. Hence we can and we will change one strictly
positive vector to another one considering asymptotic behavior if it
is necessary.

Considering maps $f_{\mset}$ we follow as close as possible to the
ideas of W.H.M.~Zijm and use his paper \cite{zij84generalized} as a
model. Notice that we cannot use Zijm's results for \ $-f_{\mset}$,
which can be expressed using maximum, because matrices should be
nonnegative and dynamics is considered on the nonnegative cone. The
problem in transferring the results obtained for $g_{\mset}$ to
$f_{\mset}$ lies in the convexity property which $f_{\mset}$ lacks.
In particular, there is no simultaneous block-triangular
decomposition, which was extremely important in
\cite{zij84generalized,sladky:bound1,sladky:bound2}. To overcome
this difficulty we show that if the set $\mset$ satisfies the
product property then there exist matrices $B$ and $C$ in $\mset$
which give the lowest and the greatest asymptotic behavior over all
matrices in $\mset$, i.e. $B^nv\preceq A^nv\preceq C^nv$ for all
$A\in\mset$ and every strictly positive vector $v$. These matrices
we call respectively \textit{$\preceq$-minimal} and
\textit{$\preceq$-maximal} for the set $\mset$. Using these notions
we study spectral properties of $f_{\mset}$. In particular, we prove
that $f_{\mset}$ possesses a strictly positive eigenvector if and
only if some (every) $\preceq$-minimal matrix possesses a strictly
positive eigenvector. The main result shows existence of nonnegative
generalized eigenvectors of $f_{\mset}$. Finally as a corollary we
describe the asymptotic behavior of each component of
$f_{\mset}^n(v)$ by showing that $f_{\mset}^n(v)\sim A^nv$ for some
(every) $\preceq$-minimal matrix $A$. We prove some other
propositions similar to the results of~\cite{zij84generalized},
sometimes assuming additional
conditions.\\

I would like to thank St{\'e}phane Gaubert for bringing my
attention to Zijm's articles and Volodymyr Nekrashevych for
helpful suggestions.

\section{Nonnegative matrices: definitions, notations, results}

We recall in this section some (usually well-known) definitions and
results that we need about nonnegative matrices (for the references
see \cite[Chapter 2]{berm:nonmatr}, \cite{sen73nonneg},
\cite[Chapter 1]{bapat}, \cite[Chapter XIII]{gantmacher}).

A matrix $A=(a_{ij})$ is called \textit{nonnegative}
(\textit{positive}) if $a_{ij}\geq 0$ ($a_{ij}>0$) for all indices
$i,j$. Denote by $A_i$ the $i$th row of the matrix $A$ and by $v_i$
the $i$th component of a vector $v$. A vector $v$ is called
\textit{strictly positive} if $v_i>0$ for all $i$. Unless otherwise
stated, all matrices will be square of a fixed dimension $N$.
Following~\cite{howard:dynpr,rw82growth,zij84generalized} the set
$\{1,2,\ldots,N\}$ is called the \textit{state space} and denoted by
$S$. If $S_1,S_2\subset S$ then we denote by $A|_{(S_1, S_2)}$ the
restriction of the square matrix $A$ to $S_1\times S_2$ and by
$v|_{S_1}$ the restriction of the vector $v$ to $S_1$.

The spectral radius of a matrix $A$ is denoted by $spr(A)$.

We say that state $i$ \textit{has access to state} $j$ if there
exists a nonnegative integer $n$ such that the $ij$th entry of $A^n$
is positive. Matrix $A$ is called \textit{irreducible} if any two
states have access to each other. In the other case $A$ is called
\textit{reducible}.

The following theorem states some important properties of square
nonnegative matrices.

\begin{teor}[{\cite[Chapter 2]{berm:nonmatr}}]\label{teor_main_nonneg_matr}
Let $A$ be a nonnegative matrix with spectral radius $\lambda$. Then
\begin{enumerate}
    \item[a)] $\lambda$ is an eigenvalue of $A$.
    \item[b)] There exists a nonnegative eigenvector $v$ associated with
    $\lambda$.
    \item[c)] If $Au\geq\sigma u$ for $u\gneq 0$
    then $\lambda\geq\sigma$.
\end{enumerate}
If moreover $A$ is irreducible then
\begin{enumerate}
    \item[d)] There exists a strictly positive eigenvector $v$ associated with
    $\lambda$ and any nonnegative eigenvector of $A$ is a scalar multiple of $v$.
    \item[e)] If $Au\geq \lambda u$ or $Au\leq \lambda u$ for $u\geq 0$ then $Au=\lambda u$.
    \item[f)] $(\sigma I-A)^{-1}>0$ for any $\sigma>\lambda$.
    \item[g)] $spr(A|_{(C,C)})<\lambda$ for any $C\varsubsetneq
    S$. If $A$ is reducible then $spr(A|_{(C,C)})\leq\lambda$ for any $C\varsubsetneq
    S$ and $spr(A|_{(C,C)})=\lambda$ for some $C\varsubsetneq S$.
    \item[h)] If $Au>\sigma u$ ($Au<\sigma u$) for $u\geq 0$
    then $\lambda>\sigma$ (respectively, $\lambda<\sigma$).
\end{enumerate}
\end{teor}

Iterations of a matrix heavily depend on its block-triangular
structure. We will describe this following
\cite{berm:nonmatr,sen73nonneg}.

A \emph{class} of a nonnegative matrix $A$ is a subset $C$ of the
state space $S$ such that $A|_{(C,C)}$ is irreducible and such that
$C$ cannot be enlarged without destroying the irreducibility. A
class $C$ is called \emph{basic} if $spr(A|_{(C,C)})=spr(A)$,
otherwise \emph{nonbasic} (when $spr(A|_{(C,C)})<spr(A)$). It
follows that for any matrix $A$ we have a partition of the state
space $S$ into classes, say $C_1, C_2,\ldots ,C_n$. Then, after
possibly permuting the states and renumbering the classes, $A$ can
be written in the form, sometimes called \emph{the Frobenius Normal
Form},
\[
A=
\begin{pmatrix}
  A_{(1,1)} & A_{(1,2)} & \ldots & A_{(1,n)} \\
                0 & A_{(2,2)} & \ldots & A_{(2,n)} \\
                0 &               0 &  \ddots & \vdots \\
                0 &               0 &     0 & A_{(n,n)}
\end{pmatrix}
\]
where $A_{(i,j)}$ denotes $A|_{(C_i,C_j)}$. Hence classes can be
partially ordered by accessibility relation. We say that a class $C$
has \textit{access to (from)} a class $C'$ if there is an access to
(from) some (or equivalently any) state in $C$ to some (or
equivalently any) state in $C'$. A class is called \emph{final} if
it has no access to any other class.

The \textit{spectral radius of a class} $C$ is the spectral radius
of $A|_{(C,C)}$.

%

The next proposition describes when a matrix $A$ has a strictly
positive eigenvector and, what is more important for the subject of
this paper, when $(A^nv)_i\sim (A^nv)_j$ for all indices $i,j$ and
any strictly positive vector $v$.

\begin{prop}\label{prop_matrix_strictly_pos_vect} Let $A$ be a
nonnegative matrix with spectral radius $\lambda$. Then the
following conditions are equivalent:
\begin{enumerate}
    \item[a)] The matrix $A$ has a strictly positive
    eigenvector.
    \item[b)] The basic classes of $A$ are precisely its final classes.
    \item[c)] $(A^nv)_i\sim \lambda^n$ for all $i$ and for some (every) vector
    $v>0$.
    \item[d)] $(A^nv)_i\sim (A^nv)_j$ for all $i,j$ and for some (every) vector
    $v>0$.
\end{enumerate}
\end{prop}
\begin{proof}
The proof of equivalence $a\textit{)}$ and $b\textit{)}$ can be
found in \cite[Theorem~3.10]{berm:nonmatr}. The proof of the rest
will follow directly from Corollary~\ref{col_matr_growth}.
\end{proof}

Also notice that if a nonnegative matrix possesses a strictly
positive eigenvector then it is associated with the spectral radius
of this matrix.

Already the last proposition indicates importance of the position of
basic and nonbasic classes of a square nonnegative matrix $A$ for
existence of a strictly positive eigenvector and behavior of its
iterations. These positions can be defined precisely by introducing
the concept of a chain. A \emph{chain of classes} of $A$ is an
ordered collection of classes $\{C_1,C_2,\ldots, C_n\}$ such that
$C_i$ has access to $C_{i+1}$, $i=1,\ldots, n-1$. The \emph{length
of a chain} is the number of basic classes it contains. The
\emph{depth of a class} $C$ of $A$ is the length of the longest
chain that starts with $C$. The \emph{degree} $\nu(A)$ of $A$ is the
length of its longest chain. Let $S_i$ be the union of all classes
of depth $i$. The partition $\{S_0, S_1, \ldots , S_{\nu}\}$ of the
state space $S$ is called the \textit{principal partition of $S$
with respect to $A$}. Principal partitions play a fundamental role
in this paper. The next result is then straight forward.

\begin{prop}\label{prop_matr_part_forgrowth}
Let $\{S_0,S_1,\ldots,S_{\nu}\}$ be the principal partition of $S$
with respect to $A$. Then, after possibly permuting the states, $A$
can be written in the form
\[
A=
\begin{pmatrix}
  A_{(\nu,\nu)} & A_{(\nu,\nu-1)} & \ldots & A_{(\nu,0)} \\
                0 & A_{(\nu-1,\nu-1)} & \ldots & A_{(\nu-1,0)} \\
                0 &               0 &  \ddots & \vdots \\
                0 &               0 &     0 & A_{(0,0)}
\end{pmatrix},
\]
where $A_{(i,j)}$ denotes $A|_{(S_i,S_j)}$. We have that
$spr(A_{(0,0)})<spr(A)$ (if $S_0$ is not empty);
$spr(A_{(i,i)})=spr(A)$ and the final classes and basic classes of
$A_{(i,i)}$ coincide for $i=1,\ldots,\nu$. Each state in $S_{i+1}$
has access to some state in $S_{i}$ for $i\geq 1$ (here
$A_{(i,i-1)}$ is non-zero).
\end{prop}

Notice that it follows from
Proposition~\ref{prop_matrix_strictly_pos_vect} that
$A|_{(S_i,S_i)}$ possesses a strictly positive eigenvector for every
$i=1,\ldots,\nu$.

We need the following useful lemma.

\begin{lemm}[{\cite[Lemma 2.5]{zij84generalized}}]\label{lemm_nonneg_matr_less_spr}
Let $A$ be a nonnegative matrix with spectral radius $\lambda$.
\begin{enumerate}
  \item[a)] If $Av\geq\sigma v$ for some real number $\sigma$ and a real vector
$v$ with at least one positive component, then $\lambda\geq\sigma$.
  \item[b)] If $Av\geq\lambda v$ with $v>0$ then every final class of $A$ is
basic and $(Av)_i=\lambda v_i$ for every $i$ in a final class of
$A$.
\end{enumerate}
\end{lemm}

Matrices which possess strictly positive eigenvectors have the
following additional properties.

\begin{lemm}\label{lemm_positive_nonfinal_class}
Let $A$ be a nonnegative matrix with spectral radius $\lambda$ which
has a strictly positive eigenvector. Let $S_1\subset S$ be the union
of all final classes of $A$. If $Au=\lambda u$ with $u|_{S_1}>0$
then $u>0$.
\end{lemm}
\begin{proof}
Let $S_2=S\setminus S_1$. Then, after possibly permuting the states,
$A$ can be written in the form:
\[
A=
\begin{pmatrix}
  A_{(S_2,S_2)} & A_{(S_2,S_1)} \\
          0     & A_{(S_1,S_1)}
\end{pmatrix}.
\]

Each class $C$ in $S_2$ has access to some state in $S_1$,
$spr(A|_{(C,C)})<\lambda$ by Proposition
\ref{prop_matrix_strictly_pos_vect}, and $\left(\lambda
I-A|_{(C,C)}\right)^{-1}>0$ by Theorem~\ref{teor_main_nonneg_matr}
item $f\textit{)}$. It follows that $\left(\lambda
I-A|_{(S_2,S_2)}\right)^{-1}A_{(S_2,S_1)}$ has a positive element in
each row. Then
\[
u|_{S_2}=\bigl(\lambda
I-A|_{(S_2,S_2)}\bigr)^{-1}A|_{(S_2,S_1)}u|_{S_1}>0.
\]

\end{proof}

\begin{lemm}[{\cite[Lemma 2.3]{zij84generalized}}]\label{lemm_dual matr}
Let $A$ be a nonnegative matrix with spectral radius $\lambda$ which
possesses a strictly positive eigenvector. Then:
\begin{enumerate}
    \item[a)] There exists a nonnegative matrix $A^{*}$ defined
    by:
    \[
    A^{*}=\lim_{n\rightarrow\infty}
    \frac{1}{n+1}\sum_{i=0}^{n}\lambda^{-i}A^i.
    \]
    We have $AA^{*}=A^{*}A=\lambda A^{*}$ and $(A^{*})^{2}=A^{*}$.
    Moreover, $a^{*}_{ij}>0$ if and only if $j$ belong to a basic
    class of $A$ and $i$ has access to $j$ under $A$.
    \item[b)] The matrix $\lambda I-A+A^{*}$ is nonsingular.
    \item[c)] If $A^{*}v=0$ for some vector $v\geq 0$ (or $v\leq 0$),
    then $v_i=0$ for every state $i$ belonging to a basic class of
    $A$.
    \item[d)] If $Av\geq \lambda v$ (or $Av\leq \lambda v$) for some
    vector $v$ then $A^{*}v\geq v$ (respectively, $A^{*}v\leq v$).
\end{enumerate}
\end{lemm}

Notice that if $A$ is a (reducible) stochastic matrix then $A^{*}$
is a limiting transition probability matrix and the inverse of
$(I-A+A^{*})$ is the so-called fundamental matrix of the respective
Markov chain.

Asymptotic behavior of iterations of a nonnegative matrix can be
studied through its generalized eigenvectors corresponding to its
spectral radius. Let $A$ be a nonnegative matrix with spectral
radius $\lambda$. The \emph{index} $\eta(A)$ of $A$ with respect to
$\lambda$ is the smallest integer $n$ such that the null spaces of
$(A-\lambda I)^n$ and $(A-\lambda I)^{n+1}$ coincide. The elements
of $\textrm{Null}(A-\lambda I)^i\setminus \textrm{Null}(A-\lambda
I)^{i-1}$ are called the \emph{generalized eigenvectors} of order
$i$. It was proved in~\cite[Theorem~3.1]{rothbl:alg}  that
$\eta(A)=\nu(A)$. Moreover, it was shown that generalized
eigenvectors can be chosen nonnegative with special strictly
positive components. More precisely (see also
\cite{handbook:rothblum,zij84generalized}):

\begin{teor}[{\cite[Theorem~3.1]{rothbl:alg}}]\label{teor_Rothblum_generalized}
Let $A$ be a nonnegative matrix with spectral radius $\lambda$. Let
$\{S_0,S_1,\ldots,S_{\nu}\}$ be the principal partition of $S$ with
respect to $A$. Then there exists a set of nonnegative generalized
eigenvectors $v^{(1)},v^{(2)},\ldots,v^{(\nu)}$ such that
\begin{eqnarray*}
Av^{(\nu)}&=& \lambda v^{(\nu)},\\
Av^{(i)} &=& \lambda v^{(i)}+v^{(i+1)}, \quad i=\nu-1,\ldots, 2,1.
\end{eqnarray*}
Moreover
\[
v^{(i)}_j>0, \quad j\in\bigcup\limits_{k=i}^{\nu}S_k \quad \mbox{
and } \quad v^{(i)}_j=0, \quad j\in\bigcup\limits_{k=0}^{i-1}S_k.
\]
\end{teor}

To give estimates on the growth of $A^nv$ we need the following
lemma.

\begin{lemm}\label{lemm_asympt_series}
For any integer $k\geq 0$ and real $\lambda, \beta> 0$ we have
asymptotic relation
\begin{eqnarray*}
\sum_{i=0}^{n-1}\lambda^{n-i}i^k\beta^{i}\sim\left\{%
\begin{array}{ll}
    n^k\beta^n, & \hbox{if $\beta>\lambda$;} \\
    n^{k+1}\lambda^n, & \hbox{if $\beta=\lambda$.}
\end{array}
\right.
\end{eqnarray*}
\end{lemm}
\begin{proof}
The asymptotic relation $\sum_{i=0}^{n-1}i^k\sim n^{k+1}$ is
standard. Then for $\lambda=\beta$:
\[
\sum_{i=0}^{n-1}\lambda^{n-i}i^k\beta^{i}=\lambda^n\sum_{i=0}^{n-1}i^k\sim
n^{k+1}\lambda^n
\]
and for $\beta>\lambda$ we have inequalities:
\begin{eqnarray*}
(n-1)^k\beta^{n-1}\leq
\sum_{i=0}^{n-1}\lambda^{n-i}i^k\beta^{i}&=&\lambda^n
\sum_{i=0}^{n-1}i^k\left(\frac{\beta}{\lambda}\right)^{i}\leq
\lambda^n
(n-1)^k\sum_{i=0}^{n-1}\left(\frac{\beta}{\lambda}\right)^{i}\leq\\
&\leq& \lambda^n
n^k\frac{\left(\frac{\beta}{\lambda}\right)^n-1}{\frac{\beta}{\lambda}-1}
\leq n^k\frac{\beta^n}{\frac{\beta}{\lambda}-1}.
\end{eqnarray*}
\end{proof}

\begin{col}\label{col_matr_growth}
Let $A$ be a nonnegative matrix with spectral radius $\lambda$. Let
$\{S_0,S_1,\ldots,S_{\nu}\}$ be the principal partition of $S$ with
respect to $A$. Then
\[
(A^nv)_k\sim n^{i-1}\lambda^n, \mbox{ for } k\in S_i, \quad i\geq
1
\]
for any strictly positive vector $v$.
\end{col}
\begin{proof}
The statement follows from a general result about asymptotic
behavior of matrix powers obtained in \cite{rothbl:sens}. We sketch
the proof to use it later.

Using the identities for the generalized eigenvectors $v^{(i)}$ from
Theorem~\ref{teor_Rothblum_generalized} and the above lemma, one can
get by induction that
\[
A^nv^{(i)}\sim\lambda^n\sum_{j=0}^{\nu-i}n^jv^{(i+j)}, \mbox{ for }
i=\nu,\ldots,1\quad \Rightarrow\quad (A^nv^{(1)})_k\sim
n^{i-1}\lambda^n, \mbox{ for } k\in S_i.
\]
Since $spr(A|_{(S_0,S_0)})<\lambda$ the $S_0$th components of a
strictly positive vector $v$ does not effect the asymptotic behavior
of $A^nv|_{S\setminus S_0}$. Hence $A^nv|_{S\setminus S_0}\sim
A^nv^{(1)}|_{S\setminus S_0}$.
\end{proof}

\textbf{Remark.} Corollary~\ref{col_matr_growth} gives us an
algorithm of finding the growth of each component of $A^nv$ (see the
detailed analysis in \cite{rothbl:sens}). For indices in $S_i$ for
$i\geq 1$ it follows directly from the corollary. For $i\in S_0$ we
consider the matrix $A|_{(S_0,S_0)}$ and its principal partition and
so on. This algorithm can be also described using chains of classes
as follows. Take a state $i$ and the corresponding class $C_i$ which
contains $i$.  Let $\beta$ be the maximum of spectral radii of
classes $C$, where $C$ runs through all classes such that $C_i$ has
access to $C$. Consider all possible chains that start at $C_i$ and
for each chain count the number of classes $C$ in this chain with
spectral radius $\beta$. Let $k$ be the maximal among such numbers.
Then $(A^nv)_i\sim n^{k-1}\beta^n$ for any strictly positive
vector~$v$. In particular, if $i$ belongs to a final class of $A$
then $(A^nv)_i\sim \beta^n$, where $\beta$ is the spectral radius of
the class.


Corollary~\ref{col_matr_growth} implies that the components of
$A^nv$ are comparable with respect to the partial order $\preceq$.
The $\preceq$-minimal possible growth of $(A^nv)_i$ over all indices
$i$ is $\sim \gamma^n$, where $\gamma$ is the spectral radius of
some final class $C$. If state $i$ has access to state $j$ then
$(A^nv)_j\preceq (A^nv)_i$. So, $(A^nv)_i\sim (A^nv)_j$ for any two
states $i$ and $j$ in the same class of $A$.

\section{Existence of strictly positive eigenvector of
$f_{\mset}$}\label{section_exist_str_pos_f}

Let $\mset$ be a finite set of nonnegative matrices. In this section
we define the notion of a $\preceq$-minimal matrix and the principal
$\preceq$-minimal partition of $S$ with respect to $\mset$. Using
these notions we give necessary and sufficient conditions for
existence of a strictly positive eigenvector for the map
$f_{\mset}(v)=\min_{A\in\mset}Av$.

Note that in general we do not have the property that for every
vector $v$ there exists $A\in\mset$ such that $f_{\mset}(v)=Av$. The
following concept eliminates this difficulty (compare
with~\cite[Section~3.1]{sen73nonneg}).

\begin{defi}
Let $\mset$ be a set of nonnegative $N\times N$ matrices. We say
that $\mset$ satisfies the \emph{product property} if for each
subset $V\subseteq S$ and for each pair of matrices $A,B\in\mset$
the matrix $C$ defined by
\[
C_i:=\left\{%
\begin{array}{ll}
    A_i, & \hbox{if $i\in V$;} \\
    B_i, & \hbox{if $i\in S\setminus V$.}
\end{array}%
\right.
\]
belongs to $\mset$.
\end{defi}

If we have any finite set $\mset_0$ of nonnegative matrices then
we can close it with respect to the product property and obtain
another finite set $\mset$. We just take all possible matrices $C$
obtained as follows: the $i$th row of $C$ is the $i$th row of some
matrix from $\mset$. Then it is easy to see that
\[
\min_{A\in\mset_0}Av=\min_{A\in\mset}Av
\]
for any vector $v$. So we can extend our given set of matrices to a
bigger one, which satisfies the product property, without changing
the map $f_{\mset}$. Moreover, for every $v$ there exists
$A=A_v\in\mset$ such that $f_{\mset}(v)=Av$.

\emph{Hence we will always assume that $\mset$ possesses the
product property.}

As was mentioned in introduction, Richard Bellman in
\cite{bellman:quasi_l_eq} considered compact sets $\mset$ of
positive matrices and proved that $g_{\mset}$ has a strictly
positive eigenvector. It was generalized to a set of irreducible
matrices in~\cite{MandlSeneta}. A simple proof of this result was
obtained by W.H.M.~Zijm \cite{zij84generalized} following the
arguments in \cite[Appendix~B]{howard:risk}. His proof also works
for the maps $f_{\mset}$.

\begin{prop}\label{prop_growth_irred_matrices}
Suppose that every matrix in the set $\mset$ is irreducible. Then
$f_{\mset}$ possesses a strictly positive eigenvector associated
with $\lambda_{\mset}=\min_{A\in\mset} spr(A)$. Moreover, it is
unique up to a scalar multiple.
\end{prop}
\begin{proof}
Take any $B\in\mset$. Let $\lambda_B$ be the spectral radius of
$B$ and let $v$ be the corresponding strictly positive
eigenvector. Find $D\in\mset$ such that
\[
Dv=\min_{A\in\mset}Av
\]
with $D_i=B_i$ if $(Bv)_i\leq (Av)_i$ for all $A\in\mset$. If $D=B$
then $f_{\mset}(v) = Bv =\lambda_Bv$ and we are done. If $D\neq B$
then $Dv\lneq Bv=\lambda_Bv$ and $\lambda_D:=spr(D)<\lambda_B$ by
Theorem~\ref{teor_main_nonneg_matr} item~$e\textit{)}$. Apply the
same procedure for the matrix $D$ with its strictly positive
eigenvector $u$ associated with $\lambda_D$. Since $\mset$ is
finite, after a finite number of steps we will reach a matrix $M$
with spectral radius $\lambda$ and eigenvector $w$ such that
\[
Mw=\min_{A\in\mset} Aw=\lambda w.
\]

Since $Aw\geq\lambda w$ for every $A\in\mset$, we get
$\lambda=\min_{A\in\mset} spr(A)$ by
Theorem~\ref{teor_main_nonneg_matr} item~$e\textit{)}$.



Let $u,v>0$ be eigenvectors of $f_{\mset}$ and let
$f_{\mset}(v)=Av=\lambda v$ and $f_{\mset}(u)=Bu=\lambda u$. Then
$spr(A)=spr(B)=\lambda$, $Bv\geq f_{\mset}(v)=\lambda v$, and
$Bv=\lambda v$ by Theorem~\ref{teor_main_nonneg_matr}
item~$e\textit{)}$. Hence by Theorem~\ref{teor_main_nonneg_matr}
item~$d\textit{)}$ the eigenvector of $f_{\mset}$ associated with
$\lambda$ is unique up to a scalar multiple.
\end{proof}

The following lemma is an important result for understanding the
asymptotic behavior of $f_{\mset}^n(v)$. It will be used throughout
the paper.

\begin{lemm}\label{lemm_min_on_matrix}
There exists $B\in\mset$ such that $B^nv\preceq A^nv$ for any
$A\in\mset$ and $v>0$.
\end{lemm}
\begin{proof}
We use induction on dimension $N$. For $N=1$ the statement is
obvious. Suppose the lemma is correct for any dimension $<N$. Let
us fix $v>0$.

For each $A\in\mset$ and $i\in S$ we can find the asymptotic
behavior of $(A^nv)_i$ using Corollary~\ref{col_matr_growth}. Define
the set $\mset'\subset\mset$ of all matrices $B$ in $\mset$ for
which there exists $i\in S$ such that $(B^nv)_i\preceq (A^nv)_j$ for
all $A\in\mset$ and $j\in S$. Note that it follows from
Corollary~\ref{col_matr_growth} that $(B^nv)_i\sim\lambda^n$ for
some real $\lambda\geq 0$. For each matrix $B\in\mset'$ define
\[
S_0(B)=\{ j\in S\;|\;(B^nv)_i\sim (B^nv)_j\sim\lambda^n\}
\]
and $S_1(B)=S\setminus S_0(B)$. Suppose some state $i$ in $S_0(B)$
has access to some state $j$ in $S_1(B)$. Then $(B^nv)_j\preceq
(B^nv)_i$. Since the asymptotic behavior of $(B^nv)_i$ is
$\preceq$-minimal for $B$, we have $(B^nv)_i\sim (B^nv)_j$. Hence
$j\in S_0(B)$ and we have a contradiction. Thus no state in $S_0(B)$
has access to any state in $S_1(B)$, which means that
$B|_{(S_0(B),S_1(B))}=0$.

Observe that the spectral radius of every class of $B$ from $S_0(B)$
is not greater than $\lambda$. If a class $C$ from $S_0(B)$ is final
then it has spectral radius $\lambda$. The converse is also true: a
class $C$ from $S_0(B)$ with spectral radius $\lambda$ is final.
Really, suppose it is not final. Then it has access to a final class
from $S_0(B)$. Thus there exists a chain which start at $C$ that
contains at least two classes with spectral radii $\lambda$. So,
$(B^nv)_i\succeq n\lambda^n$ for $i$ in $C$ by
Corollary~\ref{col_matr_growth}.

Let us show that $\mset'$ contains a matrix $B$ with the biggest set
$S_0(B)$, i.e. such that $S_0(B)\supset S_0(A)$ for any
$A\in\mset'$. It is sufficient to show that for any matrices $B$ and
$D$ from $\mset'$ there exists $E\in\mset'$ such that $S_0(E)\supset
S_0(B)\cup S_0(D)$. Define $E$ as follows: $E_i=B_i$ for $i\in
S_0(B)$ and $E_i=D_i$ for $i\not\in S_0(B)$.
\[
E=\begin{pmatrix}
  D|_{(S_1(B),S_1(B))} & * \\
  0 & B|_{(S_0(B),S_0(B))}
\end{pmatrix}
\]
Then $E\in\mset'$ and $S_0(E)\supset S_0(B)$, because
$(E^nv)_i=(B^nv)_i$ for $i\in S_0(B)$. In order to prove that
$S_0(E)$ contains $S_0(D)$, it is sufficient to prove that each
class $C$ of $E$, which belong to $S_0(D)\setminus S_0(B)$  with
spectral radius $\lambda$ is final (if it is empty we are done). By
construction $E|_{(C,C)}=D|_{(C,C)}$ and $C$ belongs to some class
$C'$ of $D$ from $S_0(D)$. If $C\neq C'$ then
$spr(D|_{(C',C')})>spr(D|_{(C,C)})=spr(E|_{(C,C)})=\lambda$ by
Theorem~\ref{teor_main_nonneg_matr} item $g\textit{)}$ and we have
contradiction with $C'\subset S_0(D)$. Thus $C=C'$ and
$E|_{(C,S\setminus C)}=D|_{(C',S\setminus C')}=0$. So $C$ is final
and our claim is proved.

Choose $B\in\mset'$ to be a matrix with the biggest set $S_0(B)$.
Denote $S_0:=S_0(B)$ and $S_1:=S_1(B)$.

If $S_0=S$ then we are done -- the matrix $B$ satisfies the
condition of the lemma. Suppose that $S_1\neq\emptyset$. The set
$\mset|_{(S_1,S_1)}$ satisfies the product property and we can apply
induction to it. So there exists $D\in\mset$ such that
$(D|_{(S_1,S_1)})^nv|_{S_1}\preceq (A|_{(S_1,S_1)})^nv|_{S_1}$ for
any $A\in\mset$. Define a matrix $E$ in the same way as above:
$E_i=B_i$ for $i\in S_0$ and $E_i=D_i$ for $i\not\in S_0$. We want
to show that it satisfies the condition of the lemma.

Again $(E^nv)_i=(B^nv)_i$ for $i\in S_0$. So  $(E^nv)_i\preceq
(A^nv)_i$ for any matrix $A\in\mset$ for $i\in S_0$. We need to
prove the previous inequality for $i\in S_1$.
\begin{eqnarray*}
(E^nv)|_{S_1}&\leq &
(D|_{(S_1,S_1)})^nv|_{S_1}+\sum_{l=1}^{n-1}(D|_{(S_1,S_1)})^{n-l}D|_{(S_1,S_0)}
(B|_{(S_0,S_0)})^lv|_{S_0}\preceq\\
&\preceq
&(D|_{(S_1,S_1)})^nv|_{S_1}+\sum_{l=1}^{n-1}(D|_{(S_1,S_1)})^{n-l}\lambda^{l}v|_{S_1}=
\sum_{l=0}^{n-1}(D|_{(S_1,S_1)})^{n-l}\lambda^{l}v|_{S_1}.
\end{eqnarray*}
Fix $i$ in $S_1$ and let $((D|_{(S_1,S_1)})^{n}v|_{S_1})_i\sim
n^k\beta^n$. Suppose $\beta<\lambda$. Then there exists $j\in S_1$
such that $((D|_{(S_1,S_1)})^{n}v|_{S_1})_j\sim \beta^n$. Then:
\[
(E^nv)_j\preceq\sum_{l=0}^{n-1}(D|_{(S_1,S_1)})^{n-l}\lambda^{l}v|_{S_1}\sim
\sum_{l=0}^{n-1}\lambda^{l}\beta^{n-l}\sim\lambda^n
\]
and therefore $j$ must be in $S_0$. We get a contradiction, hence
$\beta\geq\lambda$.

If $\beta>\lambda$ then
\begin{eqnarray*}
(E^nv)_i&\preceq&\sum_{l=0}^{n-1}(D|_{(S_1,S_1)})^{n-l}\lambda^{l}v|_{S_1}\sim
\sum_{l=0}^{n-1}\lambda^{l}(n-l)^k\beta^{n-l}\sim
n^k\beta^n\sim((D|_{(S_1,S_1)})^{n}v|_{S_1})_i\preceq\\
&\preceq &((A|_{(S_1,S_1)})^nv|_{S_1})_i\preceq (A^nv)_i
\end{eqnarray*}
for every $A\in\mset$.

Now suppose that $\beta=\lambda$. Let $C_i$ be the class of
$D|_{(S_1,S_1)}$ that contains $i$. Then $\lambda$ is the maximum
of spectral radii of $D|_{(S_1,S_1)}|_{(C,C)}$, where $C$ runs
through all classes of $D|_{(S_1,S_1)}$ such that $C_i$ has access
to $C$. Also the maximal number of classes $C$ with
$spr(D|_{(S_1,S_1)}|_{(C,C)})=\lambda$ in chains that start at
$C_i$ is $k$. If the maximum of spectral radii of $D|_{(C,C)}$,
where $C_i$ has access to $C$, is greater than $\lambda$, then
$(D^nv)_i\succeq n^{k+1}\lambda^n$. If not then the maximal number
of classes $C$ of $B$ with $spr(D|_{(C,C)})=\lambda$ in a chain
that starts at $C_i$ is at least $k+1$, otherwise there exists a
state $j$ in $S_1$ with $(E^nv)_j\sim\lambda^n$. Thus,
$(D^nv)_i\succeq n^{k+1}\lambda^n$. Notice that the above
statement is true for any matrix $A\in\mset$, i.e., if
$((A|_{(S_1,S_1)})^{n}v|_{S_1})_i\sim n^k\lambda^n$ then
$(A^nv)_i\succeq n^{k+1}\lambda^n$. Then
\begin{equation*}
(E^nv)_i\preceq\sum_{l=0}^{n-1}(D|_{(S_1,S_1)})^{n-l}\lambda^{l}v|_{S_1}\sim
\sum_{l=0}^{n-1}\lambda^{l}(n-l)^k\lambda^{n-l}\sim
n^{k+1}\lambda^n\preceq (A^nv)_i
\end{equation*}
for any $A\in\mset$. So $(E^nv)_i\preceq (A^nv)_i$ for all $i$ and
$A\in\mset$.
\end{proof}

By similar arguments one can show that there exists $C\in\mset$ such
that $A^nv\preceq C^nv$ for every $A\in\mset$.

\begin{defi}
A matrix $B\in\mset$ which satisfies Lemma \ref{lemm_min_on_matrix}
will be called \emph{$\preceq$-minimal for} $\mset$.
\end{defi}

There is a simple (but not effective) algorithm to find all
$\preceq$-minimal matrices. We find the asymptotic behavior of
$(A^nv)_i$ for all $A\in\mset$ by Corollary \ref{col_matr_growth}
and take matrices with $\preceq$-minimal growth (such matrices exist
by Lemma~\ref{lemm_min_on_matrix}).

If the spectral radius of a matrix $A$ is zero, then $A$ is
nilpotent and asymptotic behavior of $A^nv$ is trivial.
Consideration of such matrices is elementary but does not fit
precisely in the discussion below. To avoid these unnecessary
complications and without loss of generality in the sequel all
considered matrices have spectral radius $>0$.

The principal partitions, spectral radii and degrees of every two
$\preceq$-minimal matrices coincide, which follows from
Corollary~\ref{col_matr_growth} and from the fact that $A^nv\sim
B^nv$ for any $\preceq$-minimal matrices $A,B\in\mset$. Notice that
the spectral radius of a $\preceq$-minimal matrix is equal to
$\min_{A\in\mset} spr(A)$. Denote $\lambda=spr(B)$ and $\nu=\nu(B)$
for a $\preceq$-minimal matrix $B\in\mset$.

\begin{defi}
The principal partition $\{S_0,S_1,\ldots,S_{\nu}\}$ of a
$\preceq$-minimal matrix is called the \emph{principal
$\preceq$-minimal partition} of $S$ with respect to $\mset$.
\end{defi}

The following proposition gives a sufficient condition for existence
of a strictly positive eigenvector for the map $f_{\mset}$.
Moreover, it will follow from Corollary~\ref{col charact strictl
posit } that this condition is also necessary, and what is more
important for the subject of this paper that it is equivalent to the
property that all components of the iterations $f_{\mset}^n(v)$ have
the same growth.


\begin{prop}\label{prop_growth_matrix_str_posit_vect}
Suppose that some $\preceq$-minimal matrix possesses a strictly
positive eigenvector. Then $f_{\mset}$ possesses a strictly positive
eigenvector associated with $\lambda$.
\end{prop}
\begin{proof}
Note that by Proposition~\ref{prop_matrix_strictly_pos_vect} if one
$\preceq$-minimal matrix possesses a strictly positive eigenvector
then all $\preceq$-minimal matrices do.

Let $B$ be a $\preceq$-minimal matrix with strictly positive
eigenvector $v$. Apply the same procedure as in the proof of
Proposition \ref{prop_growth_irred_matrices}. Find $D\in\mset$ such
that
\[
Dv=\min_{A\in\mset}Av
\]
with $D_i=B_i$ if $(Bv)_i\leq (Av)_i$ for all $A\in\mset$. Then
$Dv\leq \lambda v$ and $D^nv\leq \lambda^nv=B^nv$. Thus $D$ is
$\preceq$-minimal, has strictly positive eigenvector and
$spr(D)=\lambda$. Since each final class of $D$ is basic,
$(Dv)_i=(\lambda v)_i$ for all $i$ in final classes by
Theorem~\ref{teor_main_nonneg_matr} item $e\textit{)}$. Hence
$D_i=B_i$ for $i$ in the final classes of $D$ and the set of final
classes of $B$ contains the set of final classes of $D$.

By Proposition~\ref{prop_matrix_strictly_pos_vect} each nonfinal
class of $D$ is nonbasic. Let $S_1\subset S$ be the union of all
final classes and let $S_2=S\setminus S_1$. Then, after possibly
permuting the states
\[
D=\begin{pmatrix}
  D|_{(S_2,S_2)} & E \\
  0 & B|_{(S_1,S_1)}
\end{pmatrix}
\]
with $spr(D|_{(S_1,S_1)})=\lambda$ and
$spr(D|_{(S_2,S_2)})<\lambda$. Define
\begin{equation}\label{eqn_prop_min_pos_eig}
u|_{S_1}=v|_{S_1} \quad\mbox{and} \quad u|_{S_2}=(\lambda
I-D|_{(S_2,S_2)})^{-1}Ev|_{S_2}.
\end{equation}
Then $Du=\lambda u$ and thus $u>0$ by Lemma
\ref{lemm_positive_nonfinal_class}. Suppose $u_i>v_i$ for some
$i\in S$. Then it follows from $Du=\lambda u$ and $Dv\leq\lambda
v$ that
\[
D|_{(S_2,S_2)}\Bigl[u|_{S_2}-v|_{S_2}\Bigr]\geq
\lambda\Bigl[u|_{S_2}-v|_{S_2}\Bigr].
\]
This contradicts $spr(D|_{(S_2,S_2)})<\lambda$ by
Lemma~\ref{lemm_nonneg_matr_less_spr}. Hence $v\geq u>0$.

By construction $u=v$ if and only if $D=B$. We can apply the same
procedure to $D$ and $u$. On each step the set of final classes of
the new matrix is contained in the set of final classes of the
previous matrix and the next eigenvector coincides with the previous
one on the states from final classes of the new matrix. Since
$\mset$ is finite, after some steps all received matrices will have
the same set of final classes and all received eigenvectors are the
same on this set. Now suppose this process will never stabilize. It
means that all received eigenvectors are different. Since $\mset$ is
finite, some matrix appears in this process at least two times with
different strictly positive eigenvectors that coincide on the final
classes of this matrix. But by~\eqref{eqn_prop_min_pos_eig}
eigenvector of a matrix is uniquely defined by its coordinates from
the final classes of this matrix. We get a contradiction. Thus after
a finite number of steps we will reach a $\preceq$-minimal matrix
$M\in\mset$ with strictly positive eigenvector $w$ such that
\[
Mw=\min_{A\in\mset}Aw=\lambda w.
\]
\end{proof}

\begin{col}
Under the conditions of
Propositions~\ref{prop_growth_irred_matrices}
or~\ref{prop_growth_matrix_str_posit_vect} the asymptotic relation
\[
(f_{\mset}^n(u))_i\sim\lambda^n
\]
holds for any strictly positive vector $u$ and $i\in S$.
\end{col}

So, if there exists a $\preceq$-minimal matrix with strictly
positive eigenvector, then the growth exponent of each component of
$f_{\mset}^n(v)$ is equal to the spectral radius of this
$\preceq$-minimal matrix.

The next proposition with $\nu=1$ gives the basis of induction for
Lemma~\ref{lemm_main2}.

\begin{prop}\label{prop_main_for_vk_on_Dk}
Let $\{S_0,S_1,\ldots, S_{\nu}\}$ be the principal $\preceq$-minimal
partition of the state space $S$ with respect to $\mset$. Then there
exists a nonnegative vector $w$ with $w|_{S_{\nu}}>0$ such that
\[
\min_{A\in\mset} Aw=\lambda w.
\]
\end{prop}
\begin{proof}
Let $B$ be any $\preceq$-minimal matrix. Since
$\{S_0,S_1,\ldots,S_{\nu}\}$ is the principal partition of $B$, the
matrix $B|_{(S_{\nu},S_{\nu})}$ possesses a strictly positive
eigenvector $v$ associated with $\lambda$. The set
$\mset|_{S_{\nu}}=\{A|_{(S_{\nu},S_{\nu})}, A\in\mset\}$ also
satisfies the product property and $B|_{S_{\nu}}$ is
$\preceq$-minimal for it. We can apply
Proposition~\ref{prop_growth_matrix_str_posit_vect} for
$\mset|_{S_{\nu}}$. There exists a strictly positive vector $u$
defined on $S_{\nu}$ such that
\[
\min_{A\in\mset} A|_{(S_{\nu},S_{\nu})}u=\lambda u.
\]
Take $w$ such that $w|_{S_{\nu}}=u$ and $w|_{S\setminus
S_{\nu}}=0$. Then $w$ satisfies the condition of the proposition.
\end{proof}

It was shown in~\cite{zij84generalized} and~\cite{sladky:bound2}
that a stronger result holds for $g_{\mset}$, which proves existence
of a simultaneous block-triangular representation of the matrices in
$\mset$ and allows one to define the ``principal partition'' of $S$
with respect to $\mset$. This partition plays a fundamental role in
those papers. This result doesn't hold for $f_{\mset}$.

\section{Generalized eigenvectors of $f_{\mset}$}

We prove in this section two lemmata from which the main result
follows immediately. The first lemma proves existence of a
solution of a set of ``nested'' functional equations. As it was
noticed in~\cite{zij84generalized}, it can be viewed as a
generalization of the Howard's policy iteration procedure
\cite{howard:dynpr}.

Let $t$ be an integer greater than 1. Suppose that for each
$A\in\mset$ we have a sequence of vectors $r_i(A)$, $i=1,\ldots,
t-1$.

\begin{lemm}\label{lemm_main_set_funct_equat}
Assume that the set of rectangular matrices
\[
\left\{ (A,r_1(A),r_2(A),\ldots,r_{t-1}(A))\;|\;A\in\mset\right\}
\]
satisfies the product property. Suppose that there exists (for all)
a $\preceq$-minimal matrix $B\in\mset$ with a strictly positive
eigenvector $v$. Suppose furthermore $B^{*}r_{t-1}(B)>0$ for any
$\preceq$-minimal matrix $B$ (here $B^{*}$ is defined in
Lemma~\ref{lemm_dual matr}). Then there exists a solution
$\{v^{(1)},\ldots,v^{(t)}\}$ of the set of functional equations:
\begin{eqnarray*}
\min_{A\in\mset}{Av^{(t)}}&=&\lambda v^{(t)}\\
\min_{A\in\mset_i}{\left\{
Av^{(i-1)}+r_{i-1}(A)\right\}}&=&\lambda v^{(i-1)}+v^{(i)}, \quad
i=2,\ldots, t,
\end{eqnarray*}
where $\mset_i$ is defined recursively by
\begin{eqnarray*}
\mset_t&:=&\{A\;|\;A\in\mset, Av^{(t)}=\lambda v^{(t)}\},\\
\mset_i&:=&\{A\;|\;A\in\mset_{i+1}, Av^{(i)}+r_{i}(A)=\lambda
v^{(i)}+v^{(i+1)}\}, \quad i=2,\ldots, t-1.
\end{eqnarray*}
Furthermore $v^{(t)}>0$.
\end{lemm}
\begin{proof}
The set of equations
\begin{eqnarray}\label{eqn_Bx_t_lemm_main}
Bv^{(t)}&=&\lambda v^{(t)},\nonumber\\
Bv^{(i)}+r_i(B)&=&\lambda v^{(i)}+v^{(i+1)}, \quad i=1,\ldots, t-1,\\
B^{*}v^{(1)}&=&0\nonumber
\end{eqnarray}
has a unique solution
\begin{eqnarray*}
v^{(t)}&=&B^{*}r_{t-1}(B),\\
v^{(i)}&=&(\lambda I-B+B^{*})^{-1}[r_{i}(B)+B^{*}r_{i-1}(B)-v^{(i+1)}], \quad i=2,\ldots, t-1,\\
v^{(1)}&=&(\lambda I-B+B^{*})^{-1}[r_1(B)-v^{(2)}].
\end{eqnarray*}
Moreover $v^{(t)}>0$. Since we have the ``extended'' product
property, there exists a matrix $D\in\mset$ such that
\begin{eqnarray*}
Dv^{(t)}&=&\min_{A\in\mset}{Av^{(t)}},\\
Dv^{(i)}+r_{i}(D)&=&\min_{A\in\mathscr{H}_{i+1}}{\left\{Av^{(i)}+r_{i}(A)\right\}},
\ i=1,\ldots, t-1,
\end{eqnarray*}
where $\mathscr{H}_i\subset\mset$ denotes the set of matrices
which minimize the right hand side of $i$th equation above. We
choose $D=B$ if $B$ satisfies above equations, i.e. if
$B\in\mathscr{H}_1$.

Then $Dv^{(t)}\leq \lambda v^{(t)}$ and thus $D$ is
$\preceq$-minimal and possesses a strictly positive eigenvector. As
above, the set of equations (\ref{eqn_Bx_t_lemm_main}) with the
matrix $D$ instead of $B$ has a unique solution
$\{u^{(1)},\ldots,u^{(t)}\}$ with $u^{(t)}>0$ and so on. We want to
show that this process will eventually stop. It is easy to see that
if $v^{(i)}$ and $u^{(i)}$ satisfy the following properties
\begin{enumerate}
    \item[(a)] $u^{(t)}\leq v^{(t)}$;
    \item[(b)] if $u^{(i)}=v^{(i)}$ for $i=k+1,\ldots ,t$ then $u^{(k)}\leq
    v^{(k)}$;
    \item[(c)] if $u^{(i)}=v^{(i)}$ for all $i=1,\ldots ,t$ then
    $D=B$,
\end{enumerate}
then, since $\mset$ is finite, after a finite number of steps we
will reach a matrix which stays intact under application of this
process. The corresponding solution of (\ref{eqn_Bx_t_lemm_main})
will satisfy the conditions of the lemma.

Let us prove (a), (b) and (c). Let $C\subset S$ be the union of
all final classes of $D$.

(a)  Using Lemma~\ref{lemm_dual matr} and construction of
$u^{(i)}$ and $v^{(i)}$ several times we get
\begin{eqnarray*}
u^{(t)}&=&D^{*}u^{(t)}=D^{*}[Du^{(t-1)}-\lambda
u^{(t-1)}+r_{t-1}(D)]= D^{*}r_{t-1}(D)\leq\\
&\leq& D^{*}[\lambda v^{(t-1)}+v^{(t)}-Dv^{(t-1)}] =
D^{*}v^{(t)}\leq v^{(t)}.
\end{eqnarray*}

(b) Now suppose $u^{(i)}\leq v^{(i)}$ for $i=k+1,\ldots, t$.
Define vectors $\psi^{(i)}$, $i=1,\ldots, t$, such that:
\begin{eqnarray*}
Dv^{(t)}&=&\lambda v^{(t)}+\psi^{(t)},\\
Dv^{(i)}+r_{i}(D)&=&\lambda v^{(i)}+v^{(i+1)}+\psi^{(i)}.
\end{eqnarray*}
From (\ref{eqn_Bx_t_lemm_main}) for the matrix $D$ and the previous
equations we get:
\begin{eqnarray}\label{eqn_B1[x0-x1]=}
D[v^{(i)}-u^{(i)}]=\lambda[v^{(i)}-u^{(i)}]+[v^{(i+1)}-u^{(i+1)}]+\psi^{(i)}.
\end{eqnarray}
Thus, $\psi^{(i)}=0$ and $Dv^{(i)}+r_{i}(D)=Bv^{(i)}+r_{i}(B)$ for
$i=k+1,\ldots, t$. Hence $B\in\mathscr{H}_{k+1}$. It follows that
$\psi^{(k)}\leq 0$ and
\begin{eqnarray}\label{eqn_B1[xk0-xk1]=}
D[v^{(k)}-u^{(k)}]&=&\lambda[v^{(k)}-u^{(k)}]+\psi^{(k)}\qquad
\Rightarrow\qquad
(\mbox{ applying $D^{*}$ })\\
D^{*}\psi^{(k)}&=&0 \nonumber.
\end{eqnarray}
Hence $\psi^{(k)}_i=0$ for $i\in C$ by Lemma~\ref{lemm_dual matr}
item $c\textit{)}$.

Consider the case $k\geq 2$. Then $\psi^{(k-1)}_i\leq 0$ for $i\in
C$ and hence $D^{*}\psi^{(k-1)}\leq 0$ by Lemma~\ref{lemm_dual matr}
item $a\textit{)}$. Applying $D^{*}$ to $(k-1)$st equation of
(\ref{eqn_B1[x0-x1]=}) we obtain:
\[
0=D^{*}[v^{(k)}-u^{(k)}]+D^{*}\psi^{(k-1)}, \quad\mbox{ but }
\quad D[v^{(k)}-u^{(k)}]\leq\lambda[v^{(k)}-u^{(k)}].
\]
Hence $[v^{(k)}-u^{(k)}]\geq
D^{*}[v^{(k)}-u^{(k)}]=-D^{*}\psi^{(k-1)}\geq 0$, because
$\psi^{(k-1)}_i\leq 0$ for $i\in C$.

For $k=1$ we have $B\in\mathscr{H}_{2}$ and since $\psi^{(1)}_i=0$
for $i\in C$ we may choose $D_i=B_i$ for $i\in C$. In this case
$D^{*}_i=B^{*}_i$ and $u^{(1)}_i=v^{(1)}_i=0$ for $i\in C$. Thus
$D^{*}v^{(1)}=0$. It follows from (\ref{eqn_B1[xk0-xk1]=}) that
\[
[v^{(1)}-u^{(1)}]\geq D^{*}[v^{(1)}-u^{(1)}]= 0.
\]

(c) As above, $\psi^{(i)}=0$ for all $i$ and hence
$B\in\mathscr{H}_{1}$. Thus $D=B$ by construction.
\end{proof}

\begin{lemm}\label{lemm_main2}
Let $\{S_0,S_1,\ldots,S_{\nu}\}$ be the principal $\preceq$-minimal
partition with respect to $\mset$. There exists a set of nonnegative
vectors $v^{(1)},v^{(2)},\ldots,v^{(\nu)}$ such that
\begin{eqnarray}\label{eqn_minA=lambda}
\min_{A\in\mset} Av^{(\nu)}&=& \lambda v^{(\nu)}, \\
\min_{A\in\mset_{i+1}} Av^{(i)} &=& \lambda v^{(i)}+v^{(i+1)},
\quad i=\nu-1,\ldots, 2,1; \nonumber
\end{eqnarray}
where
\begin{eqnarray*}
\mset_\nu&:=&\{A\;|\;A\in\mset,  Av^{(\nu)}=\lambda v^{(\nu)}\},\\
\mset_i&:=&\{A\;|\;A\in\mset_{i+1}, Av^{(i)}=\lambda
v^{(i)}+v^{(i+1)}\}, \quad i=2,\ldots, \nu-1.
\end{eqnarray*}
Moreover
\begin{equation}\label{eqn_lem_nonnegati_constr}
v^{(i)}_j>0, \quad j\in\bigcup\limits_{k=i}^{\nu}S_k \quad \mbox{
and } \quad v^{(i)}_j=0, \quad j\in\bigcup\limits_{k=0}^{i-1}S_k.
\end{equation}
\end{lemm}
\begin{proof}
By induction on $\nu$. For $\nu=1$ the result follows from
Proposition~\ref{prop_main_for_vk_on_Dk}. Suppose that the lemma
holds for $\nu<t$ and let now $\nu=t$.

Notice that
\[
\mset_{t}=\{A\;|\;A\in\mset,  Av^{(t)}=\lambda v^{(t)} \mbox{ and
} A|_{(S\setminus S_{t},S_{t})}=0 \}
\]
for any given $v^{(t)}$ such that $v^{(t)}|_{S\setminus S_{t}}=0$.
Define the set of matrices
\[
\mathscr{H}=\{A|_{(S\setminus S_t, S\setminus S_t)},
A\in\mset_{t}\}.
\]
Clearly $\mathscr{H}$ also satisfies the product property and
$B|_{(S\setminus S_{\nu},S\setminus S_{\nu})}$ is a
$\preceq$-minimal matrix for $\mathscr{H}$ for any $\preceq$-minimal
matrix $B$ for $\mset$. Thus $S_0,S_1,\ldots ,S_{\nu-1}$ is the
principal $\preceq$-minimal partition of $\mathscr{H}$. By the
induction hypothesis there exist nonnegative vectors
$u^{(1)},u^{(2)},\ldots,u^{(t-1)}$ defined on $S\setminus S_{\nu}$
such that $u^{(t-1)}_i>0$ for $i\in S_{t-1}$ and
\begin{eqnarray*}
\min_{A\in\mathscr{H}} Au^{(t-1)}&=& \lambda u^{(t-1)},\\
\min_{A\in\mathscr{H}_{i+1}} Au^{(i)} &=& \lambda
u^{(i)}+u^{(i+1)}, \quad i=t-2,\ldots, 2,1.
\end{eqnarray*}
Now we need to find vectors $v^{(1)},v^{(2)},\ldots,v^{(t)}$ such
that (\ref{eqn_minA=lambda}) holds. Let us take
\begin{eqnarray*}
v^{(i)}_j=u^{(i)}_j \quad \mbox{ and } \quad v^{(t)}_j=0
\quad\mbox{ for } j\in S\setminus S_t.
\end{eqnarray*}
Then $\mset_i\subset\{A\;|\;A\in\mset_t, A|_{(S\setminus
S_t,S\setminus S_t)}\in\mathscr{H}_{i}\}$ for $i=1,\ldots, t-1$, and
the vectors $v^{(i)}$, independent of their coordinates on $S_{t}$,
satisfy (\ref{eqn_minA=lambda}) for states in $S\setminus S_{t}$. It
remains to determine $v^{(i)}_j$ for $j\in S_t$, $i=1,\ldots, t$.
The conditions on $v^{(i)}|_{S_t}$ are the following:
\begin{eqnarray*}
\min_{A\in\mset} A|_{(S_t,S_t)}v^{(t)}|_{S_t}&=& \lambda v^{(t)}|_{S_t},\\
\min_{A\in\mset_{i+1}}
\left\{A|_{(S_t,S_t)}v^{(i)}|_{S_t}+\sum_{j=i}^{t-1}A|_{(S_t,S_j)}w^{(i)}|_{S_j}\right\}
&=& \lambda v^{(i)}|_{S_t}+v^{(i+1)}|_{S_t}, \quad i=t-1,\ldots,
2,1.
\end{eqnarray*}
Since $\{S_0,S_1,\ldots,S_{t}\}$ is the principal partition of any
$\preceq$-minimal matrix $B\in\mset$, the matrix $B|_{(S_t,S_t)}$
possesses a strictly positive eigenvector associated with $\lambda$.
Moreover $u^{(t-1)}|_{S_{t-1}}>0$. Each final class of
$B|_{(S_t,S_t)}$ has access to some state in
$B|_{(S_{t-1},S_{t-1})}$. Thus
\[(B|_{(S_t,S_{t-1})}u^{(t-1)}|_{S_{t-1}})_i>0\] for some $i$ in
every final class of $B|_{(S_t,S_t)}$. Then
$B|_{(S_t,S_t)}^{*}B|_{(S_t,S_{t-1})}u^{(t-1)}|_{S_{t-1}}>0$ for any
$\preceq$-minimal $B$ by Lemma~\ref{lemm_dual matr} item
$a\textit{)}$. We can now apply
Lemma~\ref{lemm_main_set_funct_equat} and find $v^{(i)}|_{S_t}$.

It may happened that $v^{(i)}$ does not satisfy the nonnegativity
constrains (\ref{eqn_lem_nonnegati_constr}) on $S_t$ (they satisfy
it on $S\setminus S_t$ by induction). In this case consider
\begin{eqnarray}\label{eqn_lemm_m2_u_v}
w^{(t)}&=&v^{(t)},\\
w^{(i)}&=&v^{(i)}+\alpha v^{(i+1)}, \quad i=1,\ldots, t-1.
\nonumber
\end{eqnarray}
They also satisfy (\ref{eqn_minA=lambda}) and we can choose
$\alpha$ large enough so that $w^{(i)}_j>0$ for all $j\in S_t$,
$i=1,\ldots, t$.
\end{proof}

Now we are ready to prove the main result.

\begin{teor}\label{teor_main}
Let $\{S_0,S_1,\ldots,S_{\nu}\}$ be the principal $\preceq$-minimal
partition with respect to $\mset$. Then there exists a set of
nonnegative vectors $v^{(1)},v^{(2)},\ldots,v^{(\nu)}$ such that
\begin{eqnarray*}
\min_{A\in\mset} Av^{(\nu)}&=& \lambda v^{(\nu)},\\
\min_{A\in\mset} Av^{(i)} &=& \lambda v^{(i)}+v^{(i+1)}, \quad
i=\nu-1,\ldots, 2,1.
\end{eqnarray*}
Moreover
\[
v^{(i)}_j>0, \quad j\in\bigcup\limits_{k=i}^{\nu}S_k \quad \mbox{
and } \quad v^{(i)}_j=0, \quad j\in\bigcup\limits_{k=0}^{i-1}S_k.
\]
\end{teor}
\begin{proof}
Use Lemma~\ref{lemm_main2} to find solutions
$v^{(1)},v^{(2)},\ldots,v^{(\nu)}$ of the corresponding system
(\ref{eqn_minA=lambda}). Now consider the vectors
$w^{(1)},w^{(2)},\ldots,w^{(\nu)}$ from (\ref{eqn_lemm_m2_u_v}).
It is easy to see that for $\alpha$ large enough
\[
\min_{A\in\mset_{i+1}}Aw^{(i)}=\min_{A\in\mset_{i+2}}Aw^{(i)}=\ldots=\min_{A\in\mset}Aw^{(i)},
\quad i=1,\ldots, \nu.
\]
Hence for $\alpha$ large enough the vectors
$w^{(1)},w^{(2)},\ldots,w^{(\nu)}$ satisfy the conditions of the
theorem.
\end{proof}

\begin{col}\label{col_main_growth_formin}
Let $\{S_0,S_1,\ldots,S_{\nu}\}$ be the principal $\preceq$-minimal
partition with respect to $\mset$. Then
\[
(f_{\mset}^n(v))_i \sim n^{k-1} \lambda^n, \mbox{ where } i\in S_k,
\]
for any strictly positive vector $v$ and $k\geq 1$. Moreover, for
any $\preceq$-minimal matrix $B\in\mset$
\[
(f_{\mset}^n(v))_i \sim (B^nv)_i
\]
for any strictly positive vector $v$ and $i\in S$.
\end{col}
\begin{proof}
The proof of the first part is the same as for a single matrix (see
Corollary~\ref{col_matr_growth}). Thus $(f_{\mset}^n(v))_i \sim
(B^nv)_i$ for $i\not\in S_0$. We need to prove this asymptotic
relation for $i\in S_0$.

The upper bound $f_{\mset}^n(v)\preceq B^nv$ is obvious. Define
\[
\mathscr{H}=\{ A\in\mset\;|\;A|_{(S_0,S\setminus S_0)}=0 \} \quad
\mbox{ and } \quad f_{\mset}|_{S_0}=\min_{A\in\mathscr{H}}
A|_{(S_0,S_0)}.
\]
Then $B|_{(S_0,S_0)}$ is $\preceq$-minimal for $\mathscr{H}|_{S_0}$
for any matrix $B$ $\preceq$-minimal for $\mset$. Let
$\beta=spr(B|_{(S_0,S_0)})$ for $\preceq$-minimal $B$ (notice that
$\beta<\lambda$) and let $\{S'_0,S'_1,\ldots, S'_{\nu'}\}$ be the
principal $\preceq$-minimal partition of $S_0$ with respect to
$\mathscr{H}|_{S_0}$. By Theorem~\ref{teor_main} there exist
nonnegative vectors $w^{(1)},\ldots,w^{(\nu')}$ defined on $S_0$
such that
\begin{eqnarray*}
f_{\mset}|_{S_0}(w^{(\nu')})&=&\beta w^{(\nu')}\\
f_{\mset}|_{S_0}(w^{(i)})&=&\beta w^{(i)}+w^{(i+1)}, \quad
i=\nu'-1,\ldots, 2,1,
\end{eqnarray*}
and with specified nonnegative constrains. Notice that then
$(B^nv)_i\sim n^{k-1}\beta^n$ for $i\in S'_k$, $k\geq 1$.

Let $v$ be a strictly positive vector defined on $S\setminus S_0$
such that $A|_{(S\setminus S_0,S\setminus S_0)}v\geq \lambda v$
(take for example $v^{(1)}|_{S\setminus S_0}$). Define vectors
$u^{(i)}_{\alpha_i}, i=\nu',\ldots 2,1$, such that
$u^{(i)}_{\alpha_i}|_{S\setminus S_0}=\alpha_i v$ and
$u^{(i)}_{\alpha_i}|_{S_0}=w^{(i)}$. Then
\[
f_{\mset}(u^{(i)}_{\alpha_i})=\min_{A\in\mset}\begin{pmatrix}
  \alpha_i A|_{(S\setminus S_0,S\setminus S_0)}v+A|_{(S\setminus S_0,S_0)}w^{(i)} \\
  \alpha_i A|_{(S_0,S\setminus S_0)}v+A|_{(S_0,S_0)}w^{(i)}
\end{pmatrix}=\min_{A\in\mathscr{H}} Au^{(i)}_{\alpha_i}, \quad
i=\nu',\ldots 2,1,
\]
for $\alpha_i$ large enough. Moreover we can additionally choose
$\alpha_i$ such that $\alpha_i \lambda v\geq \alpha_i\beta
v+\alpha_{i+1} v$ for $i=\nu'-1,\ldots 2,1$. Then
\begin{eqnarray*}
f_{\mset}(u^{(\nu')}_{\alpha_{\nu'}})&\geq&\begin{pmatrix}
  \alpha_{\nu'}\lambda v \\
  f_{\mset}|_{S_0}(w^{(\nu')})
\end{pmatrix}=\begin{pmatrix}
  \alpha_{\nu'}\lambda v \\
  \beta w^{(\nu')}
\end{pmatrix}\geq \begin{pmatrix}
  \alpha_{\nu'}\beta v \\
  \beta w^{(\nu')}
\end{pmatrix}=\beta u^{(\nu')}_{\alpha_{\nu'}},\\
f_{\mset}(u^{(i)}_{\alpha_{i}})&\geq&\begin{pmatrix}
  \alpha_{i}\lambda v \\
  f_{\mset}|_{S_0}(w^{(i)})
\end{pmatrix}=\begin{pmatrix}
  \alpha_{i}\lambda v \\
  \beta w^{(i)}+w^{(i+1)}
\end{pmatrix}\geq\begin{pmatrix}
  \alpha_{i}\beta v+\alpha_{i+1} v \\
  \beta w^{(i)}+w^{(i+1)}
\end{pmatrix}=
\beta u^{(i)}_{\alpha_i}+u^{(i+1)}_{\alpha_{i+1}},
\end{eqnarray*}
for $i=\nu'-1,\ldots ,1$. It follows that
$(f_{\mset}^n(u^{(1)}_{\alpha_1}))_i\succeq n^{k-1}\beta^n$ for
$i\in S'_k$, $k\geq 1$, and the lower bound is proved for $i$ in
$S\setminus S_0$ and $S_0\setminus S'_0$. We can now do the same for
the states in~$S'_0$.
\end{proof}

\begin{col}\label{col charact strictl posit }
The following conditions are equivalent:
\begin{enumerate}
    \item[a)] Function $f_{\mset}$ has a strictly positive eigenvector.
    \item[b)] Some (every) $\preceq$-minimal matrix has a strictly positive eigenvector.
    \item[c)] $(f_{\mset}^n(v))_i\sim \lambda^n$ for all $i$ and for some (every) vector
    $v>0$.
    \item[d)] $(f_{\mset}^n(v))_i\sim (f_{\mset}^n(v))_j$ for all $i,j$ and for some (every) vector
    $v>0$.
\end{enumerate}
\end{col}

\bibliographystyle{plain}

\end{document}